\documentclass{llncs}
\usepackage{amssymb,amsmath}
\newcommand{\AAA}{{\mathcal A}}

\newcommand{\AP}{\AAA_P}
\newcommand{\AQ}{\AAA_Q}

\newcommand{\comment}[1]{}

\newcommand{\F}{{\mathbb F}}

\newcommand{\id}{{\mathbf{O}}}

\newcommand{\jid}{\bf id}

\begin{document}
\title{Improved Weil and Tate
pairings for elliptic and hyperelliptic curves}
\titlerunning{Improved Weil and Tate Pairings}

\author{Kirsten Eisentr\"ager\inst{1}
\thanks{The research for this paper was done while the first author was visiting Microsoft Research.
We thank S.\ Galbraith for constructive suggestions.}
\and Kristin Lauter\inst{2}
\and Peter L.~Montgomery\inst{2}}
\authorrunning{Eisentr\"ager, Lauter, Montgomery}
\institute{School of Mathematics, Institute for Advanced Study,
           Einstein Drive,
           Princeton, NJ 08540
    \email{eisentra@ias.edu} \\
\and
           Microsoft Research, 
           One Microsoft Way, Redmond, WA 98052 \\
    \email{klauter@microsoft.com, petmon@microsoft.com}
} 
\maketitle
\begin{abstract}
  We present algorithms for computing the {\it squared} Weil and Tate
  pairings on elliptic curves and the {\it squared} Tate pairing on
  hyperelliptic curves.  The squared pairings introduced in this paper
  have the advantage that our algorithms for evaluating them are
  deterministic and do not depend on a random choice of points.  Our
  algorithm to evaluate the squared Weil pairing is about 20\% more
  efficient than the standard Weil pairing.  Our algorithm for the
  squared Tate pairing on elliptic curves matches the efficiency of
  the algorithm given by Barreto, Lynn, and Scott in the case of
  arbitrary base points where their denominator cancellation technique
  does not apply.  Our algorithm for the squared Tate pairing for
  hyperelliptic curves is the first detailed implementation of the
  pairing for general hyperelliptic curves of genus 2, and saves an estimated
  30\% over the standard algorithm.

\end{abstract}
\section{Introduction}
The Weil and Tate pairings have been proposed for use in 
cryptography, including one-round 3-way key establishment,
identity-based encryption, and short signatures~\cite{Joux2002}.
For a fixed positive integer $m$, the Weil pairing $e_m$ is a bilinear
map that sends two $m$-torsion points on an elliptic curve
to an $m$th root of unity in the field.  For elliptic curves,
the Weil pairing is a quotient of two applications of the
Tate pairing, except that the Tate pairing needs an exponentiation
which the Weil pairing omits.

For cryptographic applications, the objective is a bilinear
map with a specific recipe for efficient evaluation, and no clear way
to invert.  The Weil and Tate pairings provide such tools. Each pairing
has a practical definition which involves finding functions with
prescribed zeros and poles on the curve, and evaluating those
functions at pairs of points.

For elliptic curves, Miller \cite{Miller86} gave an algorithm for the
Weil pairing.  (See also the Appendix B to \cite{BoFr2001}, for a
probabilistic implementation of Miller's algorithm which recursively
generates and evaluates the required functions based on a random
choice of points.)  For Jacobians of hyperelliptic curves, Frey and
R\"uck \cite{Frey} gave a recursive algorithm to generate the required
functions, assuming the knowledge of intermediate functions having
prescribed zeros and poles.

For elliptic curves, we present an improved algorithm for computing
the {\it squared} Weil pairing, $e_m(P,Q)^2$.  Our deterministic
algorithm does not depend on a random choice of points for evaluation
of the pairing.  Our algorithm saves about 20\% over the standard
implementation of the Weil pairing \cite{BoFr2001}.  We use this
idea to obtain an improved algorithm for computing the {\it squared}
Tate pairing for elliptic and hyperelliptic curves. The Tate pairing
is already more efficient to implement than the Weil pairing. Our new
squared Tate pairing is more efficient than Miller's algorithm for the
Tate pairing for elliptic curves, for another 20\% saving.  For
pairings on special families of elliptic curves in characteristics $2$
and $3$, some implementation improvements were given in~\cite{GaHaSo}
and~\cite{BKLS2002}.  Another deterministic algorithm was
given in~\cite{BKLS2002}. In~\cite{BLS03}, an algorithm for the
pairing on ordinary elliptic curves in arbitrary characteristic is
given.  Our squared pairing matches the efficiency of the algorithm
in~\cite{BLS03} in the case of arbitrary base points where their
denominator cancellation technique does not apply.

For hyperelliptic curves, we use Cantor's algorithm to produce the
intermediate functions assumed by Frey and R\"uck.  We define a
squared Tate pairing for hyperelliptic curves, and use the knowledge
of these intermediate functions to implement the pairing and give an
example.  Our analysis shows that using the squared Tate pairing
saves roughly 30\% over the standard Tate pairing for genus $2$
curves.  Our algorithm for the pairing on hyperelliptic curves can be
thought of as a partial generalization of the Barreto-Lynn-Scott
algorithm for elliptic curves in the sense that we give a
deterministic algorithm which is more efficient to evaluate than the
standard one.  It remains to be seen whether some denominator cancellation
can also be achieved in the hyperelliptic case by choosing base points
of a special form as was done for elliptic curves in~\cite{BLS03}. For
a special family of hyperelliptic curves, Duursma and Lee have given a
closed formula for the pairing in~\cite{DL}, but ours is the first
algorithm for the Tate pairing on general hyperelliptic curves, and we
have implemented the genus 2 case.  The squared Weil pairing or the
squared Tate pairing can be substituted for the Weil or Tate pairing
in many of the above cryptographic applications.

The paper is organized as follows.  Section~\ref{ECW} provides
background on the Weil pairing for elliptic curves and gives the
algorithm for computing the squared Weil pairing.  Section~\ref{ECT}
does the same for the squared Tate pairing for elliptic curves.
Section~\ref{HC} presents the squared Tate pairing for hyperelliptic
curves and shows how to implement it.  Section~\ref{Ex} gives an
example of the hyperelliptic pairing.


\section{Weil pairings for elliptic curves}\label{ECW}

\subsection{Definition of the Weil pairing}
\label{ECW1}

Let $E$ be an elliptic curve over a finite
field~$\F_q$.  In the following
$\id$ denotes the point at infinity on $E$.  If $P$ is a point on $E$,
then $x(P)$ and $y(P)$ denote the rational functions mapping $P$ to
its affine $x$- and $y$-coordinates.

Let $m$ be a positive integer. We will use the Weil pairing
$e_m(\cdot,\cdot)$ definition in~\cite[p.~107]{Silverman}.  
To compute $e_m(P,Q)$, given two
distinct $m$-torsion points $P$ and $Q$ on $E$ over an extension field, 
pick two divisors
$\mathcal{A}_P$ and $\mathcal{A}_Q$ which are equivalent to $(P) -
(\id)$ and $(Q) - (\id)$, respectively, and such that $\mathcal{A}_P$
and $\mathcal{A}_Q$ have disjoint support.  Let $f_{\mathcal{A}_P}$ be
a function on~$E$ whose divisor of zeros and poles is
$(f_{\mathcal{A}_P})=m \cdot \mathcal{A}_P$.  Similarly, let
$f_{\mathcal{A}_Q}$ be a function on~$E$ whose divisor of zeros and
poles is $(f_{\mathcal{A}_Q})=m \cdot \mathcal{A}_Q$. Then
$$e_m(P,Q)=\frac{f_{\mathcal{A}_P}(\mathcal{A}_Q)}{f_{\mathcal{A}_Q}(\mathcal{A}_P)}.$$

\subsection{Rational functions needed in the evaluation of the pairing}
\label{ECW2}
Fix an integer~$m > 0$ and an $m$-torsion point~$P$ on an elliptic
curve~$E$.  Let $\AP$ be a divisor equivalent to $(P) - (\id)$.  For a
positive integer~$j$, let $f_{j,\AP}$ be a rational function on~$E$
with divisor
\begin{equation}\nonumber
(f_{j,\AP}) = j\AP - (jP) + (\id)
\end{equation}
This means that $f_{j,\AP}$ has $j$-fold zeros and poles at the
points in $\AP$, as well as a  
simple pole at~$jP$ and a simple zero at~$\mathbf{O}$, and no
other zeros or poles. 
Since $mP = \id$, it follows that $f_{m,\AP}$ has divisor $m\AP$, so in fact 
$f_{\AP} = f_{m,\AP}$.
Throughout the paper the notation $f_{j,P}$ will be used to denote
the function $f_{j,\AP}$ with $\AP = (P) - (\id)$. 

Silverman \cite[Cor. 3.5, p.\ 67]{Silverman} shows that these
functions exist.  Each $f_{i,\AP}$ is unique up to a nonzero
multiplicative scalar.  Miller's algorithm gives an iterative
construction of these functions (see for example~\cite{BKLS2002}).
The construction of $f_{1,\AP}$ depends on $\AP$.  Given $f_{i,\AP}$
and $f_{j,\AP}$, one constructs $f_{i+j,\AP}$ as the product
\begin{equation}\label{MillerRecurrence}
f_{i+j,\AP} = f_{i,\AP} \cdot f_{j,\AP} 
      \cdot \frac{g_{iP,jP}}{g_{(i+j)P}}.
\end{equation}
Here the notation $g_{U,V}$ (two subscripts) denotes the line passing
through the points $U$ and $V$ on $E$.  The notation $g_U$ (one
subscript) denotes the vertical line through $U$ and $-U$.
For more details on efficiently computing $f_{m,\AP}$, see \cite{ELM}.
\subsection{Squared Weil pairing for elliptic curves}

The purpose of this section is to construct a new pairing, which we
call the `squared Weil pairing', and which has the advantage of being
more efficient to compute than Miller's algorithm for the original
Weil pairing.  Our algorithm also has the advantage that it is
guaranteed to output the correct answer and does not depend on
inputting a randomly chosen point.  In contrast Miller's algorithm may
restart, since the randomly chosen point can cause the algorithm to
fail.

\subsection{Algorithm for $e_m(P, Q)^2$}

Fix a positive integer $m$ and the curve $E$.
Given two $m$-torsion points
$P$ and $Q$ on $E$, we want to compute $e_m(P, Q)^2$.
Start with an addition-subtraction chain for $m$. 
That is, after an initial $1$, every element in the 
chain is a sum or difference of
two earlier elements, until an $m$ appears. 
Well-known techniques give a chain of length $O(\log(m))$.
For each $j$ in the addition-subtraction chain, form a tuple
$t_j = [j P,\, jQ, \,n_j, \, d_j]$
such that
\begin{equation}\label{njSWP}
\frac{n_j}{d_j}=\frac{f_{j,P}(Q)~f_{j,Q}(-P)}{f_{j,P}(-Q)~f_{j,Q}(P)}.
\end{equation}
\noindent Start with $t_1 = [P, \, Q, \,1,\, 1]$.  
Given $t_j$ and $t_k$, this procedure gets $t_{j+k}$:

\begin{enumerate}
\item Form the elliptic curve sums $j P + k P = (j+k) P$ and 
$jQ + kQ = (j+k)Q$.
\item Find coefficients of the line 
$g_{jP,kP}(X)= c_0 + c_1 x(X) + c_2 y(X)$.
\item Find coefficients of the line
 $g_{jQ,kQ}(X) = c_0' + c_1' x(X) + c_2' y(X)$.

\item  Set
\[ \begin{split} n_{j+k} = n_j   n_k 
  (c_0 + c_1  x(Q) + c_2  y(Q))~ 
  (c_0' + c_1' x(P) - c_2' y(P)) \\
 d_{j+k} = d_j  d_k   
  (c_0 + c_1 x(Q) - c_2 y(Q))~
  (c_0' + c_1' x(P) + c_2' y(P)). \end{split}\]
\end{enumerate}
A similar construction gives $t_{j-k}$ from $t_j$ and $t_k$.  The
vertical lines through $(j+k) P$ and $(j+k) Q$ do not appear in the
formulae for $n_{j+k}$ and $d_{j+k}$, because the contributions from
$Q$ and $-Q$ (or from $P$ and $-P$) are equal.
When $j + k = m$, this simplifies to $n_{j + k} =
n_j n_k$ and $d_{j+k} = d_j d_k$, since $c_2$ and
$c_2'$ will be zero.

When $n_m$ and $d_m$ are nonzero, then the computation 
\[
\frac{n_m}{d_m}=\frac{f_{m,P}(Q)~f_{m,Q}(-P)}{f_{m,P}(-Q)~f_{m,Q}(P)}.
\]
has been successful, and we have the correct output. If, however,
$n_m$ or $d_m$ is zero, then some factor such as $c_0 + c_1 x(Q) + c_2
y(Q)$ must have vanished.  That line was chosen to pass through $j P$,
$k P$, and $(-j-k) P$, for some $j$ and $k$.  It does not vanish at
any other point on the elliptic curve.  Therefore this factor can
vanish only if $Q = j P$ or $Q = k P$ or $Q = (-j-k) P$.  In all of
these cases $Q$ will be a multiple of $P$, ensuring $e_m(P, Q) = 1$.


\subsection{Correctness proof}

\begin{theorem}[Squared Weil Pairing Formula]\label{SWPF}
Let $m$ be a positive integer.  Suppose $P$ and $Q$ are $m$-torsion points
on $E$, with neither being the identity and $P$ not equal to $\pm
Q$.  Then the squared Weil pairing satisfies
\[\frac{f_{m,P}(Q)\cdot f_{m,Q}(-P)}{f_{m,P}(-Q) \cdot 
f_{m,Q}(P)}= (-1)^m e_m(P, Q)^2.\]
\end{theorem}
\begin{proof}
Let $R_1,R_2$ be points on $E$ such that the divisors
$\mathcal{A}_P := (P+R_1)-(R_1)$ and
$\mathcal{A}_Q:= (Q+R_2)-(R_2)$ have disjoint support. Let
$\mathcal{A}_{-Q}:= (-Q+R_2)-(R_2)$.  Let $f_{\mathcal{A}_P}$ and
$f_{\mathcal{A}_Q}$ be as above.  Then
\[e_m(P,Q)=\frac{f_{\mathcal{A}_P}((Q+R_2)-(R_2))}{f_{\mathcal{A}_Q}
((P+R_1)-(R_1))}=
\frac{f_{\mathcal{A}_P}(Q+R_2)}{f_{\mathcal{A}_P}(R_2)}\cdot
\frac{f_{\mathcal{A}_Q}(R_1)}{f_{\mathcal{A}_Q}(P+R_1)}.\]
Let $g(X) = f_{m,P}(X - R_1)$. Then
$(g) = m(P + R_1) - m(R_1) = m \mathcal{A}_P = (f_{\mathcal{A}_P}),$
This implies $g(X)/f_{\mathcal{A}_P}(X)$ is constant and
$$
  \frac{f_{\mathcal{A}_P}(Q + R_2)}
       {f_{\mathcal{A}_P}(R_2)}
= \frac{g(Q + R_2)}
       {g(R_2)}
= \frac{f_{m,P}(Q + R_2 - R_1)}
        {f_{m,P}(R_2 - R_1)}.
$$
Similarly
\begin{equation} \nonumber
    \frac{f_{\mathcal{A}_Q}(R_1)}
         {f_{\mathcal{A}_Q}(P+R_1)}
  = \frac{f_{m,Q}(R_1-R_2)}
         {f_{m,Q}(P+R_1-R_2)}.
\end{equation}
Plugging these into Miller's formula gives
\begin{equation}\nonumber
e_m(P,Q) =
    \frac{f_{m,P}(Q+R_2-R_1)}
         {f_{m,P}(R_2-R_1)}
 ~  \frac{f_{m,Q}(R_1-R_2)}
         {f_{m,Q}(P+R_1-R_2)} .
\end{equation}
Using the same argument for $e_m(P,-Q)$ we obtain
\begin{gather}\nonumber
e_m(P,-Q) =
  \frac{f_{m,P}(-Q+R_2-R_1)}
       {f_{m,P}(R_2-R_1)}
~ \frac{f_{m,-Q}(R_1-R_2)}
       {f_{m,-Q}(P+R_1-R_2)}
\\ \nonumber = \frac{f_{m,P}(-Q+R_2-R_1)}
                    {f_{m,P}(R_2-R_1)}
             ~ \frac{f_{m,Q}(-R_1+R_2)}
                    {f_{m,Q}(-P-R_1+R_2)}
\end{gather}
Hence we can simplify $e_m(P,Q)^2$ to 
$$\frac{e_m(P,Q)}{e_m(P,-Q)} 
=\frac{f_{m,P}(Q+R_2-R_1) 
           ~{f_{m,Q}(R_1-R_2)}
           ~{f_{m,Q}(-P-R_1+R_2)} }
         {{f_{m,P}(-Q+R_2-R_1)
        ~ {f_{m,Q}(-(R_1-R_2))}}
        ~ {f_{m,Q}(P+R_1-R_2)}}.$$
Let $R:=R_2-R_1$.  This equation becomes
\begin{equation} \label{Rfunction}
e_m(P,Q)^2 = 
\frac{f_{m,P}(Q+R)~f_{m,Q}(-R)~f_{m,Q}(-P+R)}
     {f_{m,P}(-Q+R)~f_{m,Q}(R)~f_{m,Q}(P-R)} .
\end{equation}

Fix two linearly independent $m$-torsion points $P$ and $Q$. 
The right side of~(\ref{Rfunction}) is a rational function of $R$; call it
$\psi=\psi(R)$.  Since $f_{m,P}$ can have zeros and poles only at $P$ and
$\mathbf{O}$, and $f_{m,Q}$ can have zeros and poles only at $Q$ and
$\mathbf{O}$, this function $\psi(R)$ can have zeros or poles only at
$R=-Q$, $Q$, $P-Q$, $P+Q$, $P$, and $\mathbf{O}$.  By looking at the factors of
$\psi$ we can check that at each of these points, the value of
$\psi(R)$ is well-defined, because the zeros and poles cancel each
other out.  Since $\psi$ is a rational function on an elliptic curve
which does not have any zeros or poles, $\psi$ must be constant.
Since for certain values of $R$, $\psi(R)= e_m(P,Q)^2$, this must be
the case for all values of $R$. Hence we may in particular choose
$R=\mathbf{O}$, or equivalently $R_1=R_2$.  So let $R_1=R_2$.  By
Lemma~\ref{Peter} below,
$$\frac{f_{m,Q}(R_1-R_2)}{f_{m,Q}(-(R_1-R_2))}=(-1)^m,$$ and by
assumption $f_{m,P}$ does not have a zero or pole at $Q$ and $f_{m,Q}$
does not have a zero or pole at $P$. Hence expression (\ref{Rfunction})
simplifies to
\begin{equation}
e_m(P,Q)^2=(-1)^m~\frac{f_{m,P}( Q)~f_{m,Q}(-P) }
                  {f_{m,P}(-Q)~f_{m,Q}( P) }.
\end{equation}
\end{proof}
\begin{lemma} \label{Peter}
  Let $f : E \rightarrow \F_q$ be a rational function on $E$ with a zero
  of order $m$ (or a pole of order $-m$) at {$\id$}.  Define $g : E
  \rightarrow \F_q$ by $g(X) = f(X)/f(-X)$.  Then $g({\id})$ is finite
  and $g({\id}) = (-1)^m$.
\end{lemma}
\begin{proof}
  The rational function $h(X) = x(X)/y(X)$ has a zero of order $1$ at
  $X ={\id}$.  The function $f_1 =
  f/h^m$ has neither a pole nor a zero at $X ={\id}$, so $f_1({\id})$
  is finite and nonzero.  We check that the rational
  function $\phi(X)= h(X)/h(-X)$ has no zeros and poles on $E$. Hence
  $\phi$ is constant.  By computing $\phi(X)$ for a finite point
  $X=(x,y)$ on $E$ with $x,y\neq 0$, we see that $\phi$ is equal to
  $-1$. Hence
$$ g(X) = \frac{f(X)}{f(-X)} 
        = \frac{h(X)^m  f_1(X)}{h(-X)^m f_1(-X)} 
        =\phi(X)^m\frac{f_1(X)}{f_1(-X)}
        =(-1)^m\frac{f_1(X)}{f_1(-X)},
$$ and $g({\id}) = (-1)^m$.
\end{proof}


\subsection{Estimated savings}
\label{Savings}

In this section we compare our algorithm for the squared Weil
pairing to Miller's algorithm for the Weil pairing.  We count
operations in the underlying finite field, counting field squarings as
field multiplications throughout. This analysis assumes that we use the 
short Weierstrass form for the elliptic curve~$E$. 

In practice, some of these arithmetic operations may be over a base
field and others over an extension field.  That issue is discussed in
more detail in~\cite{GaHaSo}.  Without knowing the precise context of
the application, we don't distinguish these, although individual costs
may differ considerably.

\subsubsection{Miller's algorithm.}  Miller's algorithm chooses two points
$R_1$, $R_2$ on $E$, and lets $\mathcal{A}_P = (P+R_1) - (R_1)$ and
$\mathcal{A}_Q = (P+R_2) - (R_2)$.  Recall that in the notation of
Section~\ref{ECW1}, $f_{\mathcal{A}_P}$ is a function whose divisor is
$ m\mathcal{A}_P$.  As in Section~\ref{ECW2}, let $f_{j,
  \mathcal{A}_P}$ be a function with divisor 
$(f_{j,\mathcal{A}_P})=j(P+R_1)-j(R_1)-(jP)+(\mathbf{O}).$
This is the function $f_j$ in the notation of
\cite[p.~611f.]{BoFr2001}.  Then $f_{m, \mathcal{A}_P}=
f_{\mathcal{A}_P}$.  As pointed out in Equation~(B.1) of
\cite[p.~612]{BoFr2001}, (\ref{MillerRecurrence}) leads to the
recurrence
\begin{eqnarray}f_{i+j,\mathcal{A}_P}(\mathcal{A}_Q)=f_{i,
\mathcal{A}_P}(\mathcal{A}_Q) \cdot f_{j,
\mathcal{A}_P}(\mathcal{A}_Q) \cdot
\frac{g_{iP,jP}(\mathcal{A}_Q)}{g_{(i+j)P}(\mathcal{A}_Q)}.
\label{newrecurrence}\end{eqnarray}

During the computations, each $f_{j,\AP}(\AQ)$ is a known field
element, unlike the unevaluated functions $f_{j,\AP}$.  Since
$\AQ$ has degree~0, the value of $f_{j,\AP}(\AQ)$ is unambiguous,
whereas $f_{j,\AP}$ is defined only up to a multiplicative scalar.

To compute the Weil pairing we need
\begin{eqnarray*}
e_m(P, Q) = \frac{f_{\mathcal{A}_P}(Q+R_2)}
                 {f_{\mathcal{A}_P}(R_2)}
           ~\frac{f_{\mathcal{A}_Q}(R_1)}
                 {f_{\mathcal{A}_Q}(P+R_1)}
     = 
           \frac{f_{m,\mathcal{A}_P}(Q+R_2)}
                 {f_{m,\mathcal{A}_P}(R_2)}
           ~\frac{f_{m,\mathcal{A}_Q}(R_1)}
                 {f_{m,\mathcal{A}_Q}(P+R_1)}.
\end{eqnarray*}
For integers $j$ in an addition-sub\-traction chain for $m$,
we will construct a tuple $t_j = [jP, \, jQ,\, n_j, \, d_j]$
where $n_j$ and $d_j$ satisfy
$$\frac{n_j}{d_j}= 
 \frac{f_{j,\mathcal{A}_P}(Q+R_2)}
                 {f_{j,\mathcal{A}_P}(R_2)}
           ~\frac{f_{j,\mathcal{A}_Q}(R_1)}
                 {f_{j,\mathcal{A}_Q}(P+R_1)}.
$$
To compute $t_{i+j}$ from $t_i$ and $t_j$, 
one uses the above recurrence~(\ref{newrecurrence})
to derive the following expression for $n_{i+j} / d_{i+j}$:
\begin{gather}\nonumber
\frac{n_{i+j}}{d_{i+j}}=
\frac{n_i}{d_i}\cdot \frac{n_j}{d_j}
      \cdot \frac{g_{iP,jP}(Q+R_2)}{g_{iP,jP}(R_2)}
      \cdot \frac{g_{(i+j)P}(R_2)}{g_{(i+j)P}(Q+R_2)}
    \\      
      \cdot \frac{g_{iQ,jQ}(R_1)}{g_{iQ,jQ}(P+R_1)}
      \cdot \frac{g_{(i+j)Q}(P+R_1)}{g_{(i+j)Q}(R_1)}\label{nijdij}.
\end{gather}
To evaluate, for example, $g_{iP, jP}(Q+R_2) / g_{iP, jP}(R_2)$, start
with the elliptic curve addition $i P + j P = (i+j)P$.  This costs 1
field division and 2 field multiplications in the generic case where
$i P$ and $j P$ have distinct $x$-coordinates and neither is $\id$.
Save the slope $\lambda$ of the line
$
    g_{iP, jP}(X) = y(X) - y(i P) - \lambda (x(X) - x(i P))
$ 
through $i P$ and $j P$.  Two field multiplications suffice to 
evaluate  $g_{iP, jP}(Q+R_2)$ and  $g_{iP, jP}(R_2)$
given $Q + R_2$ and $R_2$.  No more field multiplications or divisions
are needed to compute the numerator and denominator of
$$
\frac{g_{(i+j)P}(R_2)}
     {g_{(i+j)P}(Q+R_2)}
 = \frac{x(R_2) - x((i+j)P)}
        {x(Q+R_2) - x((i+j)P)}.
$$

Repeat this once more to evaluate the last two fractions in
(\ref{nijdij}).  Overall these evaluations cost 8 field
multiplications and 2 field divisions.  We need 10 multiplications to
multiply the six fractions, for an overall cost of 18 multiplications
and 2 divisions.
  
\subsubsection{Squared pairing.} 
The squared pairing needs $n_m / d_m$ where $n_j / d_j$ is given
by~(\ref{njSWP}).  The recurrence formula is
\begin{equation}\label{nijdijSW}
\frac{n_{i+j}}{d_{i+j}} 
  = \frac{n_i}{d_i}
   ~\frac{n_j}{d_j}
   ~\frac{g_{iP, jP}(Q)}
         {g_{iP, jP}(-Q)}
   ~\frac{g_{(i+j)P}(-Q)}
         {g_{(i+j)P}(Q)}
   ~\frac{g_{iQ, jQ}(-P)}
         {g_{iQ, jQ}(P)}
   ~\frac{g_{(i+j)Q}(P)}
         {g_{(i+j)Q}(-P)} .
\end{equation}
This time the update from $t_i = [iP, \, iQ, \, n_i, \, d_i]$ and
$t_j$ to $t_{i+j}$ needs 2 elliptic curve additions.  Each elliptic
curve addition needs 2 multiplications and 1 division in the
generic case.  We can evaluate the numerator and denominator of
$$
 \frac{g_{iP, jP}(Q)}
         {g_{iP, jP}(-Q)} 
 = \frac{y(Q) - y(iP) - \lambda(x(Q) - x(iP))}
        {y(-Q) - y(iP) - \lambda(x(-Q) - x(iP))}
        $$
with only 1 multiplication, since $x(Q) = x(-Q)$.  
        
        The fraction $g_{(i+j)P}(-Q) / g_{(i+j)P}(Q)$ simplifies to 1
        since $g_{(i+j)P}(X)$ depends only on $x(X)$, not $y(X)$.
        Overall 6 multiplications and 2 divisions suffice to
        evaluate the numerators and denominators of the six fractions
        in ~(\ref{nijdijSW}).  We multiply the four non-unit fractions
        with 6 field multiplications.

Overall, the squared Weil pairing advances from $t_i$
and $t_j$ to $t_{i+j}$ with 12 field multiplications and 2 field divisions
in the generic case, compared to 18 field multiplications and 2 field
divisions for Miller's method.  When $i = j$, each algorithm needs 2
additional field multiplications due to the elliptic curve doublings.
Estimating a division as 5 multiplications, this is roughly a
20\% savings.

\section{Squared Tate pairing for elliptic curves} \label{ECT}
\subsection{Squared Tate pairing formula}
Let $m$ be a positive integer.  Let $E$ be defined over
$\F_q$, where $m$ divides $q-1$.  
Let $E(\F_q)[m]$ denote the $m$-torsion points on $E$ over $\F_q$. 
Assume $P \in E(\F_q)[m]$,
and $Q \in E(\F_q)$, with neither being the
identity and $P$ not equal to a multiple of $Q$.  The Tate pairing
$\phi_m(P,Q)$ on $E(\F_q)[m] \times E(\F_q)/m E(\F_q)$ is defined 
in~\cite{GaHaSo} as
\[ \phi_m(P,Q) := \left(f_{\mathcal{A}_P}(\mathcal{A}_Q)\right) ^{(q-1)/m},\]
with the notation and evaluation as for the Weil pairing above.  Now
we define
\[ v_m(P,Q) := \left(\frac{f_{m,P}(Q)}{f_{m,P}(-Q)}\right) ^{(q-1)/m},\]
where $f_{m,P}$ is as above, and call $v_m$ the squared Tate pairing.
To justify this terminology, we will show below that
$v_m(P,Q)= \phi_m(P,Q)^2.$

\subsection{Algorithm for $v_m(P, Q)$}
Fix a positive integer $m$ and the curve $E$.  Given an $m$-torsion
point $P$ on $E$ and a point $Q$ on $E$, we want to compute $v_m(P,
Q)$.  As before, start with an addition-subtraction chain for $m$.
For each $j$ in the chain, form a tuple $t_j =
[jP, \, n_j, \, d_j]$ such that
\begin{equation} \label{Tatequotient}
\frac{n_j}{d_j} = \frac{f_{j,P}(Q)}{f_{j,P}(-Q)}.
\end{equation}
Start with $t_1 = [P, 1, 1]$.  Given $t_j$ and $t_k$, this procedure
gets $t_{j+k}$:
\begin{enumerate}
\item Form the elliptic curve sum $j P + k P = (j+k) P$.
\item Find the line $g_{jP,kP}(X) = c_0 + c_1 x(X) + c_2 y(X)$.
\item Set 
\[  \begin{split}
n_{j+k} = n_j \cdot n_k  \cdot (c_0 + c_1  x(Q) + c_2  y(Q)) \\
d_{j+k} = d_j \cdot  d_k  \cdot (c_0 + c_1 x(Q) - c_2  y(Q)).
\end{split}\]
\end{enumerate}
A similar construction gives $t_{j-k}$ from $t_j$ and $t_k$.  The
vertical lines through $(j+k) P$ and $(j+k) Q$ do not appear in the
formulae for $n_{j+k}$ and $d_{j+k}$, because the contributions from
$Q$ and $-Q$ are equal. When $j + k = m$, one can further simplify
this to $n_{j + k} = n_j \cdot n_k$ and $d_{j+k} = d_j \cdot d_k$,
since $c_2$ will be zero.  When $n_m$ and $d_m$ are nonzero, then the
computation of (\ref{Tatequotient}) with $j=m$
is successful, and after raising to the $(q-1)/m$ power, we have the
correct output. If some $n_m$ or $d_m$ were zero, then some factor
such as $c_0 + c_1 x(Q) + c_2 y(Q)$ must have vanished.  That line was
chosen to pass through $j P$, $k P$, and $(-j-k) P$, for some $j$ and
$k$.  It does not vanish at any other point on the elliptic curve.
Therefore this factor can vanish only if $Q = j P$ or $Q = k P$ or $Q
= (-j-k) P$ for some $j$ and $k$.  In all of these cases $Q$ would be
a multiple of $P$, contrary to our assumption.

\subsection{Correctness proof}
\begin{theorem} \label{STPF}
  Let $m$ be a positive integer.  Suppose $P \in E(\F_q)[m]$ and 
$Q \in E(\F_q)$ with neither being the identity and $P \neq \pm Q$.  Then
the squared Tate pairing is
\[\phi_m(P, Q)^2 = \left(\frac{f_{m,P}(Q)}{f_{m,P}(-Q)}\right)^{(q-1)/m}.\]
\end{theorem}
\begin{proof}
Let $R_1$ and $R_2$ be as in the proof of Theorem~\ref{SWPF}.
The proof proceeds exactly as the correctness proof for the Weil
pairing. The only difference is that the
factor of $(-1)^m$ is missing in the Tate pairing and so we have
\begin{eqnarray*}
\phi_m(P,Q)^2 = \frac{\phi_m(P,Q)}{\phi_m(P,-Q)}=
\left(\frac{f_{m,P}(Q+R_2-R_1)} {f_{m,P}(-Q+R_2-R_1)}\right)^{(q-1)/m}.
\end{eqnarray*}
By the same argument as in the proof for the Weil pairing we may
choose $R_2=R_1$, which gives us the desired formula. 
\end{proof}
\subsection{Estimated savings}
This analysis is almost identical to that for the Weil
pairing in Section~\ref{Savings}.
When analyzing Miller's algorithm for the Tate pairing, the main
difference from Section~\ref{Savings} is that the analog
of~(\ref{nijdij}) has 2 fewer fractions to evaluate and combine.  An
elliptic curve addition costs 1 division and 2 multiplications, while
2 multiplications are needed to evaluate the numerators and
denominators of the two fractions.  Then 6 multiplications are needed
to combine the numerators and denominators of the 4
fractions. Therefore each step of Miller's
algorithm performing an addition costs 1 division and 10
multiplications.

For the squared Tate pairing, the analog
of~(\ref{nijdijSW}) also has 2 fewer fractions in it.  An elliptic
curve addition costs 1 division and 2 multiplications, while only 1
multiplication is needed to evaluate the numerators and denominators
of the 2 fractions.  Then 4 multiplications are needed to combine
the numerators and denominators of the 3 non-unit fractions.
Therefore each step of the squared Tate pairing
algorithm performing an addition costs 1 division and 7
multiplications.

Overall, the squared Tate pairing advances from $t_i$ and $t_j$ to
$t_{i+j}$ with 7 field multiplications and 1 field division in the generic
case, compared to 10 field multiplications and 1 field division for
Miller's method applied to the usual Tate pairing.  When $i = j$, each
algorithm needs one additional field multiplication due to the
elliptic curve doubling.  Estimating a division as 5
multiplications, this is roughly a 20\% savings.

Comparing our squared pairing to the algorithm from~\cite{BLS03}, the
algorithms are equally efficient in the case of general base points,
where there is no cancellation of denominators in their algorithm.  In
\cite{BLS03}, the authors show that if the security multiplier is even
($k=2d$) and the $x$-coordinate of the base point $Q$ lies in a
subfield $\F_{q^d}$, then the denominators in the Tate pairing
evaluation disappear.  This makes their method more efficient, but it
is possible that adding this extra structure may weaken the system for
cryptographic use.
Also, in some situations, restricting to $k$~even may not be
desirable.  

\section{Squared Tate pairing for hyperelliptic curves}\label{HC}

Let $C$ be a hyperelliptic curve of genus~$g$ given by an affine model 
$y^2 = f(x)$ with $\deg f = 2g+1$   
over a finite field $\F_q$ not of characteristic $2$. The curve $C$
has one point at infinity, which we will denote by $P_{\infty}$.
Let $J = J(C)$ be the Jacobian of $C$.
If $P=(x,y)$ is a point on $C$, then $P'$ will denote the point $P':=(x,
-y)$.  We denote the identity element of $J$ by $\jid$.

The Riemann-Roch theorem assures that each element $D$ of $J$
contains a representative of the form $A-gP_{\infty}$, where $A$ is an
effective divisor of degree $g$. In addition, we will always work
with {\it semi-reduced} representatives, which means that if
a point $P=(x,y)$ occurs in $A$ then $P':= (x,-y)$ does not occur elsewhere in
$A$.  The effective divisor representing the identity element
$\jid$ will be $g P_{\infty}$.  For an element $D$ of $J$ and integer~$i$, a
representative for $i D$ will be $A_i - gP_{\infty}$, where $A_i$ is
effective of degree $g$ and semi-reduced.

To a representative $A_i- gP_{\infty}$ we associate two
polynomials $(a_i, \, b_i)$ which represent the divisor.  The first
polynomial, $a_i(x)$, is monic and has zeros at the
$x$-coordinates of the points in the support of the divisor $A_i$.
The second polynomial, $b_i(x)$, has degree less than $\deg(a_i(x))$, and 
the graph of $y=b_i(x)$ passes through the
finite points in the support of the divisor $A_i$.

\subsection{Definition of the Tate pairing}

Fix a positive integer $m$ and assume that $\F_q$ contains 
a primitive $m$th root of unity $\zeta_m$.  The Tate pairing,
$
\phi_m : J(\F_q)[m] \times J(\F_q)/mJ(\F_q) \rightarrow \F_q^*/{\F_q^*}^m
\cong \,\langle \zeta_m \rangle ,
$
is defined in~\cite[p.~871]{Frey} explicitly as follows. 
Let $D \in J(\F_q)[m]$ and $E \in J(\F_q)$.
Let $h_{m,D}$ be a function on $C$ whose divisor is $(h_{m,D})= mD.$ Then
\[
\phi_m(D,E):= h_{m,D}(E)^{\frac{q-1}{m}} \in \langle \zeta_m \rangle.
\]
This pairing is known to be well-defined, bilinear, and non-degenerate.
The value $h_{m,D}(E)$ is defined only up to $m$th powers,
so we raise the result to the power $\frac{q-1}{m}$ to eliminate
all $m$th powers.  Note that $E$ is a divisor on the curve $C$, not an
elliptic curve. We also assume that the support of $E$ does not contain
$P_{\infty}$ and that $E$ is prime to the $A_i$'s.  
Actually $E$ needs to be prime to only those representatives
which will be used in the addition-subtraction chain for $m$, so to
about $\log m$ divisors.  

Frey and R\"uck \cite[pp.~872-873]{Frey} show how to evaluate the Tate
pairing on the Jacobian of a curve assuming an explicit reduction
algorithm for divisors on a curve.  Cantor~\cite{Cantor} gives such an
algorithm for hyperelliptic curves when the degree of $f$ is odd.  In
Section~\ref{u_ij} below, we use Cantor's algorithm to explicitly
compute the necessary intermediate functions.  These functions will be
used to evaluate the {\it squared} Tate pairing, but they could just
as well be used to evaluate the usual Tate pairing.

\subsection{Squared Tate pairing $v_m$ for
 hyperelliptic curves} 
\begin{theorem}
Given an $m$-torsion element $D$ of $J$ and
an element $E$ of $J$, with representatives 
$D= P_1+P_2+ \cdots +P_g-gP_{\infty}$ and 
$E = Q_1+Q_2+\cdots+Q_g-gP_{\infty }$
respectively,
with $P_i$ not equal to $Q_j$ or $Q'_j$ for any $i,j$ define
\[v_m(D,E):= \left(h_{m,D}
(Q_1-Q'_1+Q_2-Q'_2+\dots+Q_g-Q'_g)\right)^{(q-1)/m}.\]
Then $v_m(D,E)= \pm \phi_m(D,E)^2$
where $\phi_m(D,E)$ is the Tate pairing defined above.
\end{theorem}
\begin{proof}
  Recall that if $P_1=(x,\, y)$ is a point on $C$, then $P_1'$ is the
  point $(x,\, -y)$.
  Similarly, if $D=P_1+P_2+\dots+P_g-gP_{\infty}$, let $D' = P_1' +
  P_2'+\dots+P_g - gP_{\infty}$.  For the proof, we will compute
  $\phi_m(2D,\, 2E)$.
  
  Observe that $E-E' = Q_1-Q_1'+Q_2-Q_2'+\dots+Q_g-Q_g' \sim 2 E$ in
  the Jacobian of $C$, since
  $E+E'=(Q_1+Q_1'-2P_{\infty})+\dots+(Q_g+Q_g'-2P_{\infty})\sim \jid$. Let
  $h_{m,D}$ denote the rational function on $C$ with divisor
  $(h_{m,D})=mP_1+\dots+mP_g-2gmP_{\infty}$ as above.  Then the divisor of
  $h_{m,D}/{h_{m,D'}}$ has the form
  $$\left(\frac{h_{m,D}}{h_{m,D'}}\right) = mP_1-mP_1'+\dots+mP_g-mP_g',$$
  so $(h_{m,D}/{h_{m,D'}}) \sim 2mD$ in the Jacobian.  That means we
  can use $h_{m,D}/h_{m,D'}$ to compute the pairing $\phi_m(2D,2E)$.
  If $Q$ is any point on $C$, then we can see by comparing the divisors 
  of the two functions that $h_{m,D}(Q)= c \cdot h_{m,D'}(Q')$, where 
  $c$ is a constant which does not depend on $Q$.  

Hence 
\begin{eqnarray*}
\phi_m(2D,\, 2E)
     &=&\left(\frac{h_{m,D}(E-E')}
             {h_{m,D'}(E-E')}\right)^{(q-1)/m}
=    \left(\frac{h_{m,D}(E-E')}
             {h_{m,D}(E'-E)}\right)^{(q-1)/m} \\
&=&     \left({h_{m,D}(E-E')}^2 \right)^{(q-1)/m}.
\end{eqnarray*}
Since $\phi_m(2D,\,2E)=\phi_m(D,\,E)^4$, it follows that
\[\phi_m(D,\,E)^2 = \pm ({h_{m,D}(Q_1-Q_1'+\dots+Q_g-Q_g')})^{(q-1)/m}.\]
\end{proof}
\subsection{Functions needed in the evaluation of the pairings} 
Let $D$ be an $m$-torsion element of $J$. For a positive integer $j$, let
$h_{j,D}$ denote a rational function on $C$ with divisor \[(h_{j,D}) =
j A_1 - A_j - (j-1)g P_{\infty}.\] Since $D$ is an $m$-torsion
element, we have that $A_m=g P_{\infty}$, so the divisor of $h_{m,D}$
is $(h_{m,D}) = m A_1 - m \cdot g P_{\infty}$.  Each 
$h_{j,D}$ is well-defined up to a multiplicative constant.

Given positive divisors $A_i$ and $A_j$, we can use Cantor's algorithm to
find a positive divisor $A_{i+j}$ and a function $u_{i,j}$ with divisor 
equal to $$(u_{i,j}) = A_i+A_j-A_{i+j}-gP_{\infty}.$$
We construct $h_{j,D}(E)$ iteratively.  For $j=1$,
let $h_{1,D}$ be 1.  Suppose we have $A_i$, $A_j$, $h_{i,D}(E)$ and
$h_{j,D}(E)$. Let $u_{i,j}$ be the above function on $C$.
Then \[h_{i+j,D}(E)= h_{i,D}(E) \cdot h_{j,D}(E) \cdot u_{i,j}(E).\]

\subsection{Algorithm to compute $v_m(D,E)$} \label{u_ij}
Let $D$ and $E$ be as above.  Form
an addition-subtraction chain for~$m$.  For each $j$ in the chain we
need to form a tuple $t_j =[A_j,\,n_j,\,d_j]$ such that $jD$ has
representative $A_j - 2P_{\infty}$ and
\[
\frac{n_j}{d_j} = 
       \frac{h_{j,D}(Q_1)~h_{j,D}(Q_2)}
            {h_{j,D}(Q_1')~h_{j,D}(Q_2')}.
\]
Let $t_1=[A_1, \, 1, \, 1]$.  Given $t_i$ and $t_j$, let $(a_i, \,
b_i)$ and $(a_j, \, b_j)$ be the polynomials corresponding to the
divisors $A_i$ and $A_j$.  Do a composition step as in Cantor's
algorithm to obtain $(a,\, b)$ corresponding to $A_i + A_j$, without
performing the reduction step.  Let $d(x) = \gcd (a_i(x), \,a_j(x), \,
b_i(x)+b_j(x))$.  The output polynomials $a$, $b$, and $d$ depend on
$i$ and $j$, but we will omit the subscripts here for ease of
notation. If $d(x) = 1$, then $a(x) = a_i(x)a_j(x)$, and 
$b(x)$ is the polynomial with $\deg(b) < \deg(a)$ such that $y=b(x)$ 
passes through the distinct finite points in the support of $A_i$ and $A_j$.

The reduction step described in \cite[p.~99]{Cantor} then replaces
$(a,b)$ by $(\tilde a, \tilde b)$ where $\tilde a=(f-b^2)/a, \tilde b
\equiv -b \pmod{\tilde a}$ and $\deg(\tilde b) < \deg(\tilde a)$. This
reduction step is applied repeatedly until $\deg(\tilde a) \leq g$. In the
genus 2 situation, it follows from \cite[p.~99]{Cantor} that at most
one reduction step is performed.

\vspace{0.3cm}
\noindent
{\bf Case i.}  If $g=2$ and $\deg(a(x)) > 2$, a reduction step is performed. 
If we let
\begin{equation}\label{uijcomp}
v_{i,j}(P)= \frac{a(x(P))}{b(x(P))+y(P)},
\end{equation}
and 
\[
u_{i,j}(P):=v_{i,j}(P) \cdot d(x(P)),
\]
then $(u_{i,j}) = A_i + A_j - A_{i+j} - 2P_{\infty}$, and 
\begin{eqnarray*}
          \frac{u_{i,j}(P )}
               {u_{i,j}(P')}
      = \frac{a(x(P))}
               {a(x(P'))}
    \cdot \frac{b(x(P')) + y(P')}
               {b(x(P )) + y(P )}
    \cdot \frac{d(x(P ))}
               {d(x(P'))}
=     \frac{b(x(P')) + y(P')}
               {b(x(P)) + y(P )}.
\end{eqnarray*}
Let 
\begin{equation}\label{HCSPnij} \begin{split}
n_{i+j}:= n_{i} \cdot n_{j} \cdot (b + y)(Q_1') \cdot (b+y)(Q_2') \\
d_{i+j}:=d_{i} \cdot d_{j} \cdot(b + y)(Q_1) \cdot (b +y)(Q_2). \end{split}
\end{equation}
There is no contribution from $a$ in $n_{i+j}$ and $d_{i+j}$
because the contributions from $Q_i$ and $Q_i'$ are equal.  This improves
the algorithm for the Tate pairing in \cite{Frey}.

\vspace{0.3cm}
\noindent
{\bf Case ii.} If $g=2$ and $\deg(a(x)) \leq 2$, then
$u_{i,j}(P) = d(x(P))$. In this case we let $n_{i+j}:= n_{i} \cdot
n_{j}$ and $d_{i+j}:=d_{i} \cdot d_{j}$.

\vspace{0.3cm}
\noindent
{\bf Case iii.}
Suppose $g > 2$. If $r$ reduction steps are needed, then 
to compute $u_{i,j}$, we obtain intermediate factors
$v^{(1)}_{i,j},\dots,v^{(r)}_{i,j}$, 
one factor as in (\ref{uijcomp}) per
reduction step. Then $u_{i,j}$ will be the product
$u_{i,j}:= v^{(1)}_{i,j}\cdot{} \dots {}\cdot v^{(r)}_{i,j}\cdot d(x(P))$.

{\bf Note:} If we evaluate $n_i$ and $d_i$ at intermediate steps 
then it is not enough
to assume that the divisors $D$ and $E$ are coprime.  Instead, $E$
must also be coprime to $A_i$ for all $i$ which occur in the addition
chain for $m$.  One way to ensure this condition is to require that
$E$ and $D$ be linearly independent and that the polynomial $p(x)$ in
the pair $(p(x), \, q(x))$ representing $E$ be irreducible.  There are
other ways possible to achieve this, like changing the addition chain
for $m$.
\subsection{Estimated savings for genus 2}
Using a straightforward implementation of Cantor's algorithm,
the total costs for doubling and addition on the Jacobian of a
hyperelliptic curve of genus 2 in odd characteristic, $C:y^2= f(x)$,
where $f$ has degree $5$, are as follows. Doubling an element costs
$34$ multiplications and $2$ inversions. Adding two distinct elements
of $J$ costs $26$ multiplications and $2$ inversions. More efficient
implementations of the group law may alter the total impact of our
algorithm.  Different field multiplication/inversion ratios and field
sizes, as well as differing costs in an extension field will also
affect the analysis, but these costs are chosen as representative for
the purpose of estimating the savings.
 
\subsubsection{Analysis of standard algorithm}
Let $D:= P_1 + P_2 - 2 P_{\infty}$.
Let $R_1$, $R_2$, $R_3$, $R_4$ be four points on $C$ such that
$Q_1 + Q_2 - 2 P_{\infty} \sim R_1 +R_2 - R_3 - R_4$
in $J$.  The algorithm in \cite{Frey} computes $t_{i+j}$ from
$t_i$ and $t_j$, where $t_i = [A_i, \, n_j, \, d_j]$ and
$$ 
\frac{n_j}{d_j} 
        = \frac{h_{j,D}(R_1)
                 ~ h_{j,D}(R_2)}                 
                 { h_{j,D}(R_3)
                 ~ h_{j,D}(R_4)}.
$$
The expression for $n_{i+j} / d_{i+j}$ becomes
\begin{eqnarray}
\nonumber
\frac{n_{i+j}}{d_{i+j}}
    =  \frac{n_i}{d_i}
     ~ \frac{n_j}{d_j}
     ~ \frac{u_{i,j}(R_1)
           ~  u_{i,j}(R_2)}
            {u_{i,j}(R_3)
           ~  u_{i,j}(R_4)} \phantom{.}.
   & \label{nijdijHC}
\end{eqnarray}
To form $u_{i,j}$, we have to perform an addition or doubling step to
obtain $A_{i+j}$ from $A_i$ and $A_j$.  This costs 34 multiplications
and 2 inversions for a doubling, 26 multiplications and 2 inversions
for an addition.  Then
\[u_{i,j}(P)= \frac{a(x(P))}{b(x(P)) + y(P)},\] 
and to compute $(n_{i+j}, \, d_{i+j})$, we need to evaluate $u_{i,j}$
at four different points. Each evaluation of $a(x(P))$ costs 2
multiplications in a doubling step,
3 multiplications in an addition step (square or product of monic quadratics). 
Evaluation of $b(x(P))$ (cubic)
costs 3 multiplications. Finally we multiply the partial numerators
and denominators out, using 5 multiplications each, including the
multiplications with $n_i$, $n_j$, $d_i$, and $d_j$. So the total cost
for an addition step is 60 multiplications and 2 inversions, and the
total cost for a doubling is 64 multiplications and 2 inversions.
\subsubsection{Squared Tate pairing}
The squared Tate pairing works with the divisor $Q_1-Q_1'+Q_2-Q_2'
\sim 2Q_1 + 2Q_2 - 4 P_{\infty}$.  After adding $A_i$ and $A_j$ to
obtain $A_{i+j}$ as above, we need to form
\begin{eqnarray}
\nonumber
\frac{n_{i+j}}{d_{i+j}}
    =  \frac{n_i}{d_i}
     ~ \frac{n_j}{d_j}
     ~ \frac{u_{i,j}(Q_1)
           ~  u_{i,j}(Q_1')}
            {u_{i,j}(Q_2)
           ~  u_{i,j}(Q_2')} \phantom{.}.
   & \label{nijdijHCSP}
\end{eqnarray}
As can be seen from~(\ref{HCSPnij}) above, no evaluations of $a(x(P))$
are needed.  For $i=1,2$, we need to evaluate $b(x(Q_i))$ and
$b(x(Q_i'))$.  This costs only 3 multiplications for each $i$, since
the $x$-coordinates of $Q_i$ and $Q_i'$ are the same.  Finally, we
have to multiply the partial numerators and denominators, for a total
cost of 12 multiplications for either a doubling or an addition.

So the total cost for an addition step is 38 multiplications and 2
inversions, and the total cost for a doubling is 46 multiplications
and 2 inversions.  Estimating an inversion as 4 multiplications, this
is a 25\% improvement in the doubling case and a 33\% improvement in
the addition case.

\section{Example: $g=2$, $p=31$, $m=5$} \label{Ex}

In this section, we evaluate the squared Tate pairing on $5$-torsion
on the Jacobian of a hyperelliptic genus $2$ curve over a field of
$31$ elements. Let $C$ be defined by the affine model $y^2=f(x)$ where
$f(x) = x^5 +13x^4+2x^3+4x^2+11x+1.$ The group of points on the
Jacobian of $C$ over $\F_{31}$ has order $N = 1040$.  Let $D$ be the
$5$-torsion element of the Jacobian of $C$ given by the pair of
polynomials $D = [x^2+23x + 15, 13 x + 28].$ Let $E$ be the element of
the Jacobian of $C$ of order $260$ given by the pair $E=[x^2+4x+2, 29
x + 20].$ Then the squared Tate pairing evaluated at $D$ and $E$ is
$v_5(D,E)=4,$ where
$$h_{5,D}= \frac{(x + 26)^2(x^4+19x^3+23x^2+16x
  +19)(x^2+23x+15)}{x^3+6x^2+9 x + 21 + y}.$$

To illustrate the bilinearity of the pairing, look for example at
$2D = [x^2  + 25 x + 9, 10 x + 6],$
$3D = [x^2  + 25 x + 9, 21 x + 25],$ and
$2E = [x^2  + x + 3, 26 x + 3].$
Then we compute that indeed $v_5(2D,E)=16 = v_5(D,E)^2,$
with
$$h_{5,2D}= \frac{(x + 26)(x^4+19x^3+23x^2+16x
  +19)^2(x^2+25x+9)}{(x^3+6x^2+9 x + 21 + y)^2},$$
and $v_5(D,2E) =
16= v_5(D,E)^2,$ with $h_{5,D}$ as above.  Also $$v_5(3D,E)\equiv 2
\equiv v_5(D,E)^3 \pmod{31},$$
with
$$h_{5,3D}= \frac{(x + 26)(x^4+19x^3+23x^2+16x +19)^2(x^2+
  25x+9)}{(30x^3+25x^2+22 x + 10 + y)^2}.$$

\end{document}